\documentclass[12pt]{amsart}
\usepackage{amsfonts,amsmath,amsthm,amssymb}
\usepackage{url}

\theoremstyle{remark}

\theoremstyle{definition}

\newcommand{\FF}{\mathbb{F}}
\newcommand{\ZZ}{\mathbb{Z}}

\newcommand{\RR}{\mathbb{R}}

\newcommand{\cF}{\mathcal{F}}

\DeclareMathOperator{\kernel}{Ker}

\usepackage{lscape}

\begin{document}
\title{On the classification of self-dual $\ZZ_k$-codes II}

\author{Masaaki Harada}
\address[Corresponding author]{Research Center for Pure and Applied
Mathematics, 
Graduate School of Information Sciences, 
Tohoku University, Sendai 980--8579, Japan}
\email{mharada@m.tohoku.ac.jp}
\author{Akihiro Munemasa}
\address{Research Center for Pure and Applied Mathematics, 
Graduate School of Information Sciences, 
Tohoku University, Sendai 980--8579, Japan}
\email{munemasa@math.is.tohoku.ac.jp}
\date{September 24, 2015}

\begin{abstract}
In this short note, we report the classification of
self-dual $\ZZ_k$-codes of length $n$ for
$k \le 24$ and $n \le 9$.
\end{abstract}
\maketitle

\section{Introduction}\label{sec:Intro}
Let $\ZZ_{k}$ be the ring of integers modulo $k$, where $k$
is a positive integer greater than $1$.
A $\ZZ_{k}$-code $C$ of length $n$
is a $\ZZ_{k}$-submodule of $\ZZ_{k}^n$.
A code $C$ is {\em self-dual} if $C=C^\perp$, where
the dual code $C^\perp$ of $C$ is defined as 
$C^\perp = \{ x \in \ZZ_{k}^n\ | \ x \cdot y = 0$ for all $y \in C\}$
under the standard inner product $x \cdot y$. 
Two $\ZZ_k$-codes $C$ and $C'$ are {\em equivalent} 
if there exists a monomial $(\pm 1, 0)$-matrix $P$ with 
$C' = C \cdot P$, 
where $C \cdot P = \{ x P \mid x \in C\}$.  
A {\em Type~II} $\ZZ_{2k}$-code was defined in~\cite{BDHO} 
as a self-dual code with the property that all
Euclidean weights are divisible by $4k$
(see~\cite{BDHO} for the definition of Euclidean weights).
It is known that a Type~II $\ZZ_{2k}$-code of length $n$ exists 
if and only if $n$ is divisible by eight~\cite{BDHO}.
A self-dual code which is not Type~II is called {\em Type~I}.

As described in~\cite{RS-Handbook},
self-dual codes are an important class of linear codes for both
theoretical and practical reasons.
It is a fundamental problem to classify self-dual codes.
Much work has been done towards classifying self-dual $\ZZ_k$-codes
for small $k$ and  modest $n$ (see~\cite{RS-Handbook}).
Let $n_{\max}(k)$ denote the
maximum integer $n$ such that self-dual $\ZZ_k$-codes are classified
up to length $n$.
For $k =2,3,\ldots,10$, we list in Table~\ref{Tab:K}
our present state of knowledge about $n_{\max}(k)$.
We also list the reference for the classification of
self-dual $\ZZ_k$-codes of length $n_{\max}(k)$.

\begin{table}[tbh]
\caption{Known classification of self-dual $\ZZ_k$-codes} 
\label{Tab:K}
\begin{center}
{\small
\begin{tabular}{c|ccccccccc}
\noalign{\hrule height0.8pt}
$k$ & 2 & 3 & 4 & 5 & 6 & 7 & 8 & 9 & 10 \\
\hline
$n_{\max}(k)$ & 40 & 24 & 19 & 16 & 12 & 12 & 12 & 12 & 10 \\
\hline
Reference & \cite{BDM} & \cite{HM-GF3} & \cite{HM09} & \cite{HO-GF5} 
& \cite{HM09} & \cite{HO-GF7} & \cite{HM09} & \cite{HM09} &\cite{HM09} \\
\noalign{\hrule height0.8pt}
\end{tabular}
}
\end{center}
\end{table}

A classification method of self-dual
$\ZZ_k$-codes based on a classification of $k$-frames of
unimodular lattices was given by
the authors and Venkov~\cite{HMV}.
Then, in~\cite{HM09}, using this method,
self-dual $\ZZ_k$-codes were classified for $k=4,6,8,9,10$ (see Table~\ref{Tab:K}).
Using the same method, in this short note, 
we complete the classification of 
self-dual codes $\ZZ_k$-codes of length $n$ for
$k \le 24$ and $n \le 9$.
All computer calculations in this short note
were done by {\sc Magma}~\cite{Magma}.

\section{Classification of self-dual $\ZZ_k$-codes}

\subsection{Method for classifications}

A classification method of self-dual
$\ZZ_k$-codes based on a classification of $k$-frames of
unimodular lattices was given by
the authors and Venkov~\cite{HMV}.
We describe it briefly here (see~\cite{HM09} and~\cite{HMV} for 
undefined terms and details).  

 
A set $\{f_1, \ldots, f_{n}\}$ of $n$ vectors $f_1, \ldots, f_{n}$ in an
$n$-dimensional unimodular lattice $L$ with
$(f_i, f_j) = k \delta_{i,j}$
is called a {\em $k$-frame} of $L$,
where $(x,y)$ denotes the standard inner product of $\RR^n$, and 
$\delta_{i,j}$ is the Kronecker delta.
The following construction of lattices from codes
is called {\em Construction A\@}.
If $C$ is a self-dual $\ZZ_k$-code of length $n$
then
\[
A_{k}(C) = \frac{1}{\sqrt{k}}
\{(x_1,\ldots,x_n) \in \ZZ^n \:|\:
(x_1 \bmod k,\ldots,x_n \bmod k)\in C\}
\]
is an $n$-dimensional unimodular lattice.
Moreover, $C$ is Type~II if and only if
$A_k(C)$ is even.
Let $\cF=\{f_1,\dots,f_n\}$ be a $k$-frame of $L$. Consider the
mapping
\begin{align*}
&\pi_{\cF}:\frac1k\bigoplus_{i=1}^n\ZZ f_i\to\ZZ_k^n\\
&\pi_{\cF}(x)=((x,f_i) \bmod k)_{1\leq i\leq n}.
\end{align*}
Then $\kernel\pi_{\cF}=\bigoplus_{i=1}^n\ZZ f_i\subset L$, so
the code $C=\pi_{\cF}(L)$ satisfies $\pi_{\cF}^{-1}(C)=L$.
This implies $A_k(C) \simeq L$, and every code $C$ with
$A_k(C) \simeq L$ is obtained as $\pi_{\cF}(L)$ for some
$k$-frame $\cF$ of $L$, 
where $L \simeq L'$ means that $L$ and $L'$ are isomorphic lattices.
Moreover,
every Type~I (resp.\ Type~II) $\ZZ_k$-code of length $n$
can be obtained from a certain $k$-frame in some $n$-dimensional
odd (resp.\ even) unimodular lattice.

Let $L$ be an $n$-dimensional unimodular lattice, and let
$\cF=\{f_1,\dots,f_n\}$, $\cF'=\{f'_1,\dots,f'_n\}$ be
$k$-frames of $L$. Then the self-dual codes $\pi_{\cF}(L)$ and
$\pi_{\cF'}(L)$ are equivalent if and only if there exists an
automorphism $P$ of $L$ such that
$\{\pm f_1,\dots,\pm f_n\}\cdot P=
\{\pm f'_1,\dots,\pm f'_n\}$~\cite{HMV}.
This implies that the classification of codes $C$ satisfying 
$A_k(C) \simeq L$ reduces to finding a set of representatives of 
$k$-frames in $L$ up to the action of the automorphism group of $L$.


\subsection{Results}

Here, we report the classification of
self-dual $\ZZ_k$-codes of length $n$ for
$k \le 24$ and $n \le 9$.
Our classification method of self-dual $\ZZ_k$-codes of length $n$
requires a classification of
$n$-dimensional unimodular lattices.
For $n \le 7$, any $n$-dimensional unimodular lattice
is isomorphic to $\ZZ^n$.
Up to isomorphism, there are two 
$8$-dimensional unimodular lattices, one of which is the even 
unimodular lattice denoted by $E_8$ and the other is $\ZZ^8$.
Also, up to isomorphism, there are two 
$9$-dimensional unimodular lattices,
$\ZZ^9$ and $E_8 \oplus \ZZ$ (see~\cite[p.~49]{SPLAG}).

In Table~\ref{Tab:result}, we list the number of 
inequivalent self-dual
$\ZZ_k$-codes $C$ with $A_k(C) \simeq L$ for $k \in \{2,3,\ldots,24\}$
and $L \in \{\ZZ^i \mid i=1,2,\ldots,9\} \cup \{E_8,E_8\oplus\ZZ\}$.
Note that all self-dual
$\ZZ_k$-codes $C$ with $A_k(C) \simeq E_8$ are Type~II.
A classification of self-dual $\ZZ_k$-codes
of lengths $n \le 9$ was known for some $k$.
In this case, we list the references in the last columns of
the table.
Generator matrices can be obtained electronically from~\cite{Data}.
All the zero entries in Table~\ref{Tab:result} are explained as follows.
For $k \in \{3,6,7,11,12,14,15,19,21,22,23,24\}$, 
if there is a self-dual $\ZZ_k$-code of length $n$, then
$n$ is divisible by four (see~\cite[Corollary 2.2]{DHS}).
For $k \in \{2,5,8,10,13,17,18,20\}$,
if there is a self-dual $\ZZ_k$-code of length $n$,
then $n$ is even  
(see~\cite[Theorem 4.2]{DGW06}, \cite[Corollary 2.2]{DHS}).
If $k$ is a square, then 
there is a self-dual $\ZZ_k$-code for every length
(see~\cite{C-S-Z4}, \cite{DGW06}).
If a self-dual $\ZZ_k$-code is Type~II, then $k$ is even.

\begin{table}[tbh]
\caption{Classification of self-dual $\ZZ_k$-codes of lengths $n \le 9$} 
\label{Tab:result}
\begin{center}
{\footnotesize
\begin{tabular}{c|ccccccccccc|c}
\noalign{\hrule height0.8pt}
$k$ & $\ZZ$ &$\ZZ^2$ & $\ZZ^3$ & $\ZZ^4$ & $\ZZ^5$ 
& $\ZZ^6$ & $\ZZ^7$ & $\ZZ^8$ & $E_8$ & $\ZZ^9$ &$E_8\oplus\ZZ$ & Reference\\
\hline
 2 &0& 1 & 0 & 1 & 0 & 1 &  0 &  1 &  1 &  0&   0& \cite{Pless72}\\
 3 &0& 0 & 0 & 1 & 0 & 0 &  0 &  1 &  0 &  0&   0& \cite{MPS}\\
 4 &1& 1 & 1 & 2 & 2 & 3 &  4 &  7 &  4 &  7&   4& \cite{C-S-Z4,Ga96}\\
 5 &0& 1 & 0 & 1 & 0 & 2 &  0 &  3 &  0 &  0&   0& \cite{LPS-GF5}\\
 6 &0& 0 & 0 & 1 & 0 & 0 &  0 &  3 &  2 &  0&   0& \cite{DHS,HM09,KO,Park09}\\
 7 &0& 0 & 0 & 1 & 0 & 0 &  0 &  4 &  0 &  0&   0& \cite{PT}\\
 8 &0& 1 & 0 & 1 & 0 & 3 &  0 & 20 &  9 &  0&   0& \cite{DGW06,HM09}\\
 9 &1& 1 & 2 & 3 & 3 & 6 &  9 & 16 &  0 & 28&   7& \cite{BBN,HM09}\\
10 &0& 1 & 0 & 2 & 0 & 5 &  0 & 16 & 11 &  0&   0& \cite{HM09}\\
11 &0& 0 & 0 & 1 & 0 & 0 &  0 &  8 &  0 &  0&   0& \cite{BGGHK}\\
12 &0& 0 & 0 & 2 & 0 & 0 &  0 & 73 & 22 &  0&   0& \\
13 &0& 1 & 0 & 2 & 0 & 5 &  0 & 21 &  0 &  0&   0& \cite{BGGHK}\\
14 &0& 0 & 0 & 1 & 0 & 0 &  0 & 27 & 18 &  0&   0& \\
15 &0& 0 & 0 & 2 & 0 & 0 &  0 & 51 &  0 &  0&   0& \\
16 &1& 1 & 1 & 2 & 3 & 7 & 23 &295 & 63 &697& 141& \\
17 &0& 1 & 0 & 2 & 0 & 6 &  0 & 47 &  0 &  0&   0& \cite{BGGHK}\\
18 &0& 1 & 0 & 4 & 0 &12 &  0 &178 & 69 &  0&   0& \\
19 &0& 0 & 0 & 2 & 0 & 0 &  0 & 57 &  0 &  0&   0& \\
20 &0& 1 & 0 & 2 & 0 &17 &  0 &725 &176 &  0&   0& \\
21 &0& 0 & 0 & 3 & 0 & 0 &  0 &208 &  0 &  0&   0& \\
22 &0& 0 & 0 & 2 & 0 & 0 &  0 &166 & 75 &  0&   0& \\
23 &0& 0 & 0 & 1 & 0 & 0 &  0 &120 &  0 &  0&   0& \\
24 &0& 0 & 0 & 1 & 0 & 0 &  0 &3690 &456 &  0&   0\\
\noalign{\hrule height0.8pt}
\end{tabular}
}
\end{center}
\end{table}

\subsection{Remark on length 4}
A classification of self-dual $\ZZ_k$-codes
of length $4$ was given in~\cite{BGGHK} for $k=19,23$, 
and in~\cite{Park} for prime $k \le 100$.
We note that the definition of equivalence employed in~\cite{Park}
is different from our definition.
Let $N_4(k)$ denote the number of inequivalent self-dual $\ZZ_k$-codes
of length $4$.
We give in Table~\ref{Tab:result4}
the numbers $N_4(k)$ for integers $k$ with $25 \le k \le 200$.
We remark that the classification can be extended
to $k=1000$.
However, in order to save space, we do not list the result.

\begin{table}[tbhp]
\caption{Classification of self-dual $\ZZ_k$-codes of length $4$
$(25 \le k \le 200)$} 
\label{Tab:result4}
\begin{center}
{\footnotesize
\begin{tabular}{c|c||c|c||c|c||c|c||c|c||c|c}
\noalign{\hrule height0.8pt}
$k$ & $N_4(k)$ &$k$ & $N_4(k)$ &$k$ & $N_4(k)$ &$k$ & $N_4(k)$ &
$k$ & $N_4(k)$ &$k$ & $N_4(k)$ \\
\hline
 25&  5& 55&  5& 85& 10&115&  9&145& 14&175& 20\\
 26&  3& 56&  1& 86&  6&116&  5&146& 11&176&  2\\
 27&  4& 57&  7& 87&  7&117& 15&147& 18&177& 14\\
 28&  3& 58&  5& 88&  2&118&  8&148&  8&178& 13\\
 29&  2& 59&  3& 89&  5&119&  8&149&  7&179&  8\\
 30&  5& 60&  5& 90& 19&120&  5&150& 30&180& 19\\
 31&  2& 61&  4& 91&  9&121&  9&151&  7&181&  9\\
 32&  1& 62&  4& 92&  3&122&  9&152&  3&182& 19\\
 33&  4& 63&  8& 93&  8&123& 11&153& 20&183& 15\\
 34&  4& 64&  2& 94&  6&124&  6&154& 15&184&  3\\
 35&  3& 65&  8& 95&  8&125& 13&155& 12&185& 17\\
 36&  6& 66&  9& 96&  1&126& 20&156& 14&186& 20\\
 37&  3& 67&  4& 97&  6&127&  6&157&  8&187& 14\\
 38&  3& 68&  4& 98& 10&128&  1&158& 10&188&  6\\
 39&  5& 69&  5& 99& 13&129& 12&159& 12&189& 26\\
 40&  2& 70&  9&100& 12&130& 21&160&  2&190& 23\\
 41&  3& 71&  3&101&  5&131&  6&161& 10&191&  8\\
 42&  5& 72&  4&102& 14&132&  9&162& 27&192&  2\\
 43&  3& 73&  5&103&  5&133& 11&163&  8&193& 10\\
 44&  2& 74&  6&104&  3&134&  9&164&  7&194& 14\\
 45&  7& 75& 11&105& 16&135& 22&165& 25&195& 31\\
 46&  3& 76&  5&106&  8&136&  4&166& 11&196& 16\\
 47&  2& 77&  5&107&  5&137&  7&167&  7&197&  9\\
 48&  2& 78& 10&108&  9&138& 15&168&  5&198& 33\\
 49&  6& 79&  4&109&  6&139&  7&169& 15&199&  9\\
 50& 10& 80&  2&110& 14&140&  9&170& 26&200& 10\\
 51&  6& 81& 12&111& 10&141& 10&171& 21&&\\
 52&  5& 82&  7&112&  3&142&  9&172&  8&&\\
 53&  3& 83&  4&113&  6&143& 10&173&  8&&\\
 54&  8& 84&  9&114& 14&144&  6&174& 20&&\\
\noalign{\hrule height0.8pt}
\end{tabular}
}
\end{center}
\end{table}

Let $s_1,s_2,\ldots ,s_u$ be positive integers.
An orthogonal design of order $n$ and of type
$(s_1,s_2,\ldots ,s_u)$, denoted $OD(n;s_1,s_2,\ldots
,s_u)$, on the commuting variables $x_1,x_2,\ldots ,x_u$ is an 
$n \times n$ matrix $A$ with entries from $\{0, \pm x_1, \pm x_2,
\ldots,\pm x_u\}$ such that
\[
AA^T=\left( \sum_{i=1}^{u}s_ix_i^2 \right) I_n,
\]
where $A^T$ denotes the transpose of $A$ and $I_n$ is the identity
matrix of order $n$. 
The following matrix
\[
M(x_1,x_2,x_3,x_4)=
\left( \begin{array}{rrrr}
 x_1& x_2& x_3& x_4\\
-x_2& x_1&-x_4& x_3\\
-x_3& x_4& x_1&-x_2\\
-x_4&-x_3& x_2& x_1
\end{array}
\right)
\]
is well known as an $OD(4;1,1,1,1)$.
From Lagrange's theorem on sums of squares,
for each positive integer $k$,  the matrix
$M$ gives a $k$-frame of $\ZZ^4$. 
However, there are $k$-frames which are not obtained in this way.
Indeed, if $k$ is a square, then a $k$-frame can be obtained from a
$k$-frame of $\ZZ^3$, for example,
\[
{\cF}_9=
\{
( 1, 2, 2,0),
(-2,-1, 2,0),
(-2, 2,-1,0),
( 0, 0, 0,3)
\}
\]
is a $9$-frame.
Although the following matrix
\[
N(x_1,x_2,x_3,x_4)=
\left( \begin{array}{rrrr}
 x_1& x_2& x_3&x_4\\
-x_2& x_1&-x_4&x_3\\
 x_4&-x_3& x_1&x_2\\
 x_3& x_4&-x_2&x_1
\end{array}
\right)
\]
is not an orthogonal design,
if $x_1x_3+x_1x_4-x_2x_3+x_2x_4 = 0$ then 
\[
N(x_1,x_2,x_3,x_4)
N(x_1,x_2,x_3,x_4)^T
=( \sum_{i=1}^{4}x_i^2 )I_4.
\]
A $15$-frame ${\cF}_{15}$ is obtained from $N(3,1,2,-1)$.
We also found the following $21$-frame ${\cF}_{21}$:
\[
{\cF}_{21}=
\{
( 4, 1, 0, 2),
( 0,-4, 1, 2),
( 1, 0, 4,-2),
(-2, 2, 2, 3)
\}.
\]
Note that $N_4(9)=3$, $N_4(15)=2$ and $N_4(21)=3$.
The two other $9$-frames are obtained from $M(3,0,0,0)$ and $M(2,2,1,0)$.
The other $15$-frame is obtained from $M(3,2,1,1)$.
The two other $21$-frames are obtained from 
$M(0, 1, 2, 4)$ and $M(2, 2, 2, 3)$.


\subsection{Remark on length 8}

Let $N_{8,I}(2k)$ (resp.\ $N_{8,II}(2k)$) be the number of 
inequivalent Type~I
(resp.\ Type~II) $\ZZ_{2k}$-codes of length $8$.
From Table~\ref{Tab:result}, we see
$N_{8,I}(2) = N_{8,II}(2)$ and
$N_{8,I}(2k) > N_{8,II}(2k)$ $(k=2,3,\ldots,12)$.
We conjecture that $N_{8,I}(2k) > N_{8,II}(2k)$
for all integers $k$ with $k \ge 2$.

\bigskip
\noindent {\bf Acknowledgment.}
This work is supported by JSPS KAKENHI Grant Number 26610032.


\end{document}